\documentclass{alf-j-l}
\usepackage{verbatim} 
\usepackage{url}  
\input AlfPreambleL  
\newtheorem{theorem}{Theorem}

\newtheorem{proposition}[theorem]{Proposition} 

\theoremstyle{definition}

\newtheorem{Example}{Example} 
\newtheorem{example*}{Example}

\theoremstyle{remark}

\newtheorem*{remark*}{Remark}

\PII{Centre for Number Theory Research, Macquarie University, Sydney}%

\copyrightinfo{\number\year}{Alfred J van der Poorten}

\begin{document}

\def\currentvolume{165}
\def\currentissue{Draft of\ }
\def\paperdate{\today}
\def\ISSN{}
\pagespan{69}{\pageref{page:lastpage}}   

\def\ceNTRelogo{\relax}
\title {Quadratic Irrational Integers\rlap{\smash{\,\qquad\ceNTRelogo}}\\
With Partly Prescribed \\Continued Fraction
Expansion}

\author{Alfred J. van der Poorten}
\address{ceNTRe
for Number Theory Research, 1 Bimbil Place, Sydney, NSW 2071, 
\penalty10000 Australia}
\email{alf@math.mq.edu.au (Alf van der Poorten)}

\dedicatory{To the memory of Professor Dr. Bela Brindza}

\thanks{This paper was commenced while the author enjoyed the
hospitality and support of the Erwin Schr\"odinger Institute, Vienna and
substantially augmented during an immediately subsequent visit to the
Kossuth Lajos University, Debrecen. The author was supported in part by
a grant from the Australian Research Council.}


\date{\Now.}


\begin{abstract}  We generalise remarks of Euler and of
Perron by explaining how to detail all quadratic integers for which the
symmetric part of their \cf\
\ex\ commences with prescribed partial quotients.\end{abstract}

\maketitle
\pagestyle{myheadings}\markboth{{\headlinefont 
Alf van der Poorten}}{{\headlinefont Partly presribed quadratic integers}}

\noindent I last saw Bela Brindza, my once postdoctoral student,  in
April, 2002. I was working on the paper below and attempted to enthuse
him with its results, particularly those concerning periodic \ex s in
function fields of characteristic~zero.

\section{Periodic Continued Fractions}

\noindent Suppose $\omega$ satisfies $\omega^2-t\omega+n=0$ and is an integer. That is, its
trace $t=\omega+\ov\omega$ and norm $n=\omega\ov\omega$ both are rational integers. Then
it is well known that the \cf\ \ex\ of $\omega$ is periodic and is of the shape
\begin{equation} \label{eq:symmetry}
\omega=\CF{a_0\\\ov{a_1\\\ldots\\a_{r-1}\\2a_0-(\omega+\ov\omega)}},
\end{equation}
where the word $\cCF{a_1\\\ldots\\a_{r-1}}$ is a palindrome.

\begin{Example} \label{ex:1}
For example, $\sqrt{61}=\CF{7\\\ov{1\\4\\3\\1\\2\\2\\1\\3\\4\\1\\14}}$. Just so,
$(1+\sqrt{61}\,)/2$, an integer of trace~$1$, has \ex\ $\CF{4\\\ov{2\\2\\7}}$.
\end{Example}

Although the question was already asked and is partially answered by Euler,
see the very interesting translation \cite{WW}, and is discussed by
Perron
\cite{Pe}, it is noticeably less well known that given an arbitrary palindrome
$\cCF{a_1\\\ldots\\a_{r-1}}$ in positive integers there are infinitely many positive
integers $A=a_0$ so that
\eqref{eq:symmetry} displays the \ex\ of a quadratic integer. 

We explain that argument and
rather more. Indeed, noting that for every $h$ we have the \ex\
$\omega=\CF{A\\a_1\\\ldots\\a_h\\(\omega+P_{h+1})/Q_{h+1}}$ with
integers $P=P_{h+1}$, and $Q=Q_{h+1}$ (that is, the complete quotients
$\omega_{h+1}$ of $\omega$ all are of the indicated shape), we find all
quadratic integers $\omega$ for which the symmetric part of its \cf\ \ex\ commences with
the integers $\cCF{a_1\\\ldots\\a_{h}}$. Specifically, we find the constraints on $A$,
$P$, and $Q$ so that $\omega$ is indeed integral with trace
$t$ and norm $n$.

Our remarks are of particular interest in the function field case, where
such translations as `positive integer'$=$ `\poly\ of degree at least one',
`integer part' $=$ `\poly\ part', and `integer'
$=$ `\poly' are to be applied. There, however, the results in the the characteristic two
case demand a distinct summary.

\section{Continued Fractions
}
\noindent  Anyone attempting to compute the truncations
$\CF{a_0\\a_1\\\ldots\\a_h}=x_h/y_h$ of a \cf\ will be delighted to notice that
the definition
$$\CF{a_0\\a_1\\\ldots\\a_h}=a_0+ 1/\CF{a_1\\\ldots\\a_h}$$ immediately
implies by induction
on $h$ that there is a correspondence
$$\begin{pmatrix} a_0&1\\1&0\end{pmatrix}\begin{pmatrix}
a_1&1\\1&0\end{pmatrix}\cdots
\begin{pmatrix} a_h&1\\1&0\end{pmatrix}=\begin{pmatrix}
x_h&x_{h-1}\\y_h&y_{h-1}\end{pmatrix}
\longleftrightarrow \CF{a_0\\a_1\\\ldots\\a_h}=x_h/y_h
$$
between products of certain two by two matrices and the convergents of \cf s. 

If $\alpha=
\CF{a_0\\a_1\\a_2\\\ldots}$ then the sequence $(\alpha_h)$ of complete quotients
of $\alpha$ is defined by $\alpha=\CF{a_0\\a_1\\\ldots\\a_h\\\alpha_{h+1}}$.
It follows from the correspondence that
\begin{multline*}
\alpha\longleftrightarrow
\begin{pmatrix} a_0&1\\1&0\end{pmatrix} \begin{pmatrix} a_1&1\\1&0\end{pmatrix}
\cdots \begin{pmatrix} a_h&1\\1&0\end{pmatrix} \begin{pmatrix}
\alpha_{h+1}&1\\1&0\end{pmatrix} \\ =\begin{pmatrix}
x_h&x_{h-1}\\y_h&y_{h-1}\end{pmatrix} \begin{pmatrix}
\alpha_{h+1}&1\\1&0\end{pmatrix}
\longleftrightarrow \frac{x_h\alpha_{h+1}+x_{h-1}}{y_h\alpha_{h+1}+y_{h-1}}\,.
\end{multline*}
That is, we have
$$\alpha=\frac{x_h\alpha_{h+1}+x_{h-1}}{y_h\alpha_{h+1}+y_{h-1}}, \quad\text{and so}\quad 
\alpha_{h+1}=-\frac{y_{h-1}\alpha-x_{h-1}}{y_{h}\alpha-x_{h}}.
$$
Recalling that $x_{-1}=1$, $y_{-1}=0$ because an empty matrix
product is the identity matrix, we obtain
\begin{equation}(-1)^{h+1}\alpha_1\alpha_2\cdots \alpha_{h+1}=(y_{h}\alpha-x_{h})^{-1}.
\label{eq:identity}\end{equation}

\begin{Example} \label{ex:2} 
Indeed, set $\delta=\sqrt{61}$ and let $\ov\delta$ denote its conjugate, $-\sqrt{61}$. 
The \cf\
\ex\ of $\delta$ commences
\begin{alignat*}{3}
&\delta \,&&=7 &&-(\,\ov\delta+7)/1\\
\delta_1=(&\delta+7)/12&&=1 &&-(\,\ov\delta+5)/12\\
\delta_2=(&\delta+5)/3&&=4 &&-(\,\ov\delta+7)/3\\
\delta_3=(&\delta+7)/4&&=3 &&-(\,\ov\delta+5)/4\,.
\end{alignat*}
Note, in part to set notation, that a
typical line in this tableau is
\begin{equation} \label{eq:expansion}
\delta_h=(\delta+P_h)/Q_h=a_h-(\ov\delta+P_{h+1})/Q_h\,,
\end{equation} 
where 
\begin{equation} \label{eq:recursion}
P_h+P_{h+1}+(\delta+\ov\delta\,)=a_hQ_h   \quad\text{and}\quad
-Q_hQ_{h+1}=\delta\ov\delta+(\delta+\ov\delta\,)P_h+P_h^2\,.\end{equation}
Thus, by \eqref{eq:identity} and \eqref{eq:recursion},
\begin{equation} \label{eq:norm}x_h^2-61y_h^2=\Norm(x_h-\delta
y_h)=\Norm(\delta_1\delta_2\cdots\delta_{h+1})=(-1)^{h+1}Q_{h+1}\,.
\end{equation}
\end{Example}
One can, by \eqref{eq:expansion} and \eqref{eq:recursion}, readily confirm
that the
\cf\
\ex\ of
$\delta$ \emph{is} eventually periodic.  Indeed, by induction one
notices that both  
\begin{equation}\label{eq:ineq}
0\lt 2P_{h+1}+(\delta+\ov\delta)<\delta-\ov\delta \quad\text{and}\quad
0\lt Q_h<\delta-\ov\delta\,.
\end{equation}
Thus the box principle guarantees
periodicity. 

Moreover, conjugation is an involution on the tableau (hence my
eccentric way of denoting the remainders), turning the tableau upside
down. At the cost of replacing $\delta_0$ by $\delta+7$ at line $h=0$
and continuing the expansion to see the additional symmetry at lines $h=5$
and $6$ (see Example~\ref{ex:3}) we could both again have proved
periodicity and have confirmed that the preperiod is just the line $h=0$. 

Suppose, more generally, that $\omega^2-t\omega+n=0$ and that $\omega$ is
the larger real zero. Then our remarks apply without essential change to
$\omega$ in place of $\sqrt{61}$. In the function field case, with $\omega$
replaced by $Y(X)$ and $\deg_X Y=g+1$, the cited inequalities \eqref{eq:ineq}
become
$\deg Q_h\le g$ and $\deg P_{h+1}=g+1$. Now, however, the box principle does
not apply (unless the base field is finite) and periodicity is at best
happenstance. Nonetheless, if periodicity happens to happen then the period
will indeed commence no later than at line $1$. 

The issue of periodicity is encapsulated by the following remark emphasising
that periodicity coincides with the existence of units.
\begin{proposition} Let $x$ and $y$ be positive integers so that
$x^2-txy+ny^2=\pm1$. \label{pr:1}
Then the decomposition
$$N=\begin{pmatrix} x&-ny\\y&x-ty\end{pmatrix}=\begin{pmatrix} b_0&1\\1&0\end{pmatrix}
\begin{pmatrix} b_1&1\\1&0\end{pmatrix}\cdots\begin{pmatrix} b_{r}&1\\1&0\end{pmatrix}
\begin{pmatrix} 0&1\\1&0\end{pmatrix}
$$
in positive integers $b_0$, $b_1$, $\ldots\,$, $b_{r}$ entails that 
$$\omega=\CF{\ov{b_0\\b_1\\\ldots\\b_{r-1}\\b_r\\0}}=
\CF{b_0\\\ov{b_1\\\ldots\\b_{r-1}\\b_r+b_0}}
$$
for some $\omega$ satisfying $\omega^2-t\omega+n=0$.
\end{proposition}
\begin{proof}
Consider the
periodic \cf\
$$\CF{\overline{b_0\\b_1\\\ldots\\b_r\\0}\,}=\gamma\,,\quad\text{say.}
$$
Thus $\gamma=\CF{{b_0\\b_1\\\ldots\\b_r\\0}\\\gamma\,}$. By the matrix
correspondence we have
$$\gamma\longleftrightarrow N\begin{pmatrix}\gamma&1\\1&0\end{pmatrix}=
\begin{pmatrix} x\gamma-ny&x\\y\gamma+x-ty&y\end{pmatrix} \quad\text{so}\quad
\gamma=\frac{x\gamma-ny}{y\gamma+x-ty}\,.
$$
But this is $y(\gamma^2-t\gamma+n)=0$. Hence $y\ne0$ 
confirms our claim.
\end{proof}
Note here that $x/y=\CF{b_0\\b_1\\\ldots\\b_{r-1}}$ provides the unit
$x-\omega y$. Moreover, the symmetry of the matrix
$$N\begin{pmatrix} t&1\\1&0\end{pmatrix}=
\begin{pmatrix} -ny+tx&x\\x&y\end{pmatrix}
$$
confirms that the word $\cCF{b_0\\b_1\\\ldots\\b_{r-1},b_r+t}$ is a
palindrome, again explaining the symmetry under conjugation to which we
alluded above. 

\begin{Example} \label{ex:3}
The expansion of $\delta=\sqrt{61}$ continues with the lines
\begin{alignat*}{3}
\delta_4=(&\delta+5)/9&&=1 &&-(\,\ov\delta+4)/9\\
\delta_5=(&\delta+4)/5&&=2 &&-(\,\ov\delta+6)/5\\
\delta_6=(&\delta+6)/5&&=2 &&-(\,\ov\delta+4)/5
\end{alignat*}
with $Q_5=Q_6$ signalling the symmetry. Notice that in this case the
period length $r=11$ is odd (equivalently, the fundamental unit of
$\Z[\sqrt{61}]$ has norm $-1$).

For a case $r=2s$, even, the symmetry is given by the line
$$\omega_s=(\omega+P_s)/Q_s=a_s-(\ov\omega+P_{s})/Q_s\,,
$$
that is, by $P_{s+1}=P_s$\,. We have $Q_s\Div
(2P_s+t)$. In this `ambiguous' case, that is equivalent to $Q_s$ also
dividing the discriminant $t^2-4n$ of the quadratic order
$\Z[\omega]$.
\end{Example}

\section{Quadratic Integers with Partly Prescribed Expansion}
\begin{theorem}\label{th:2} Suppose that $\omega^2-t\omega+n=0$ and its \cf\
\ex\ is
$$\omega=\CF{A\\a_1\\a_2\\\ldots\\a_h\\(\omega+P_{h+1})/Q_{h+1}}.
$$
Set $P=P_{h+1}$, $Q=Q_{h+1}$, and
$$\begin{pmatrix} a_1&1\\1&0\end{pmatrix}\begin{pmatrix} a_2&1\\1&0\end{pmatrix}
\cdots\begin{pmatrix} a_{h}&1\\1&0\end{pmatrix}=
\begin{pmatrix} p&p'\\q&q'\end{pmatrix}.
$$
Then both
\begin{align*}p(P^2+tP+n)+p'Q(A+P)&=-q'Q\,,\\
p(A^2-tA+n)+q(A+P)&=-q'Q\,.
\end{align*}

\end{theorem}
\begin{proof} The data is equivalent to
\begin{multline*}
\begin{pmatrix} A&1\\1&0\end{pmatrix}\begin{pmatrix} p&p'\\q&q'\end{pmatrix}\begin{pmatrix}
1&P\\0&Q\end{pmatrix} =\begin{pmatrix} x&-ny\\y&x-ty\end{pmatrix}\\ =\begin{pmatrix}
Ap+q&Ap'+q'\\p&p'\end{pmatrix}\begin{pmatrix} 1&P\\0&Q\end{pmatrix} =\begin{pmatrix}
Ap+q&P(Ap+q)+Q(Ap'+q')\\p& Pp+Qp'\end{pmatrix}
\end{multline*} Thus $x=Ap+q$ and $y=p$  and $Pp+Qp'=Ap+q-pt$. Hence
\begin{align} \label{eq:first} p(A-t-P)-p'Q&=-q\\ p(PA+n)+p'QA&=-qP-q'Q\notag
\end{align} Strategic additions of multiples of the first equation to the
second provides the allegation.\end{proof}

\subsection{The Classical Result}\label{ss:classical}
Specifically, if $h=r-1$ so that we are given the complete symmetric
part of the period of $\omega$, then $P=A-t$, $Q=1$ and the identities
become
$$Np-Tq=-q' \quad\text{with $N=\Norm(\omega-A)$ and
$T=\Trace(\omega-A)$.}
$$
Because $pq'-p'q=(-1)^{r-1}$, we know immediately that the general solution
for the integers $N$ and $T$ is
\begin{equation}\label{eq:classical}
N=(-1)^r({\smash{q'}}^2-Lq) \quad\text{and}\quad T=(-1)^r(p'q'-Lp),
\end{equation}
where $L$ is an arbitrary integer. The latter equation
congenially reports that
$$2A=t-(-1)^r(p'q'-Lp)\,.
$$
Thus if $p'$ and $q'$ both are odd while $Lp$ is even then $t$ must be
taken odd; say $t=1$ as in $\omega=\frac12(1+\sqrt D)$; here $D=T^2-4N$.
In all other cases one may take
$t=0$ and obtain
$\omega$ the square root $\sqrt D$ of an integer, not a square; then
$4D=T^2-4N$. In either case, of course, $D=D(L)$ is a \poly\ quadratic in
the integer parameter~$L$. Distinguishing the cases brings difficulties,
see for example \cite{McL}.

In the function field case this parity matter is no issue, unless the
characteristic is $2$ where we see that then $t=p'q'+Lp$, with $L$ an
arbitrary polynomial, gives us all possible discriminants
$t^2-4n=t^2$.

\subsection{The General Result} 
We recall $pq'-p'q=(-1)^h$ and 
$(P^2+tP+n)=-Q_{h}Q:=-Q'Q$. Thus  we
have the two general solutions
\begin{alignat*}{3} P^2+tP+n &=(-1)^{h+1}Q(q'{}^2-K_1p') &\qquad &  
-(A+P) &&=(-1)^{h+1}(qq'-K_1p)\,,\\ A^2-tA+n &=(-1)^{h+1}(q'{}^2Q-K_2q)
&\qquad&  -(A+P) &&=(-1)^{h+1}(p'q'Q-K_2p)\,,
\end{alignat*} where $K_1$ and $K_2$ denote arbitrary integers. Set
$kq'=K=:K_2-K_1$.
\label{ss:general}

The dexter
pair of equations is
$Kp=q'(p'Q-q)$. Indeed $q'$ does divide $K$ and
also $q+kp\equiv0\pmod {p'}$, showing that $K$ is fixed modulo~$p'q'$.
In fact, because $pq'=(-1)^h+p'q$, we have $k\equiv(-1)^{h+1}qq'\pmod{p'}$.
Thus, for integers~$L$,
\begin{equation}\label{eq:k}
(-1)^{h+1}k=qq'-Lp' \quad\text{and then}\quad (-1)^{h+1}Q=q^2-Lp\,;
\end{equation}
and $L=q'Q-kq$. Note that \eqref{eq:k} in any case is immediate from
\eqref{eq:first} at page~\pageref{eq:first} above. Hence, at very
first glance surprisingly, knowing
$Q_{h+1}$ as well as
$a_1$,
$\ldots\,$, $a_h$ restricts us to just one parameter families of solutions.

The
sinister equations and those for
$A+P$, very helpfully yield
$A-t-P=k$, and therefore, with $T=\Trace(\omega-A)$, $N=\Norm(\omega-A)$,
\begin{equation}\label{eq:trace}
T=(-1)^{h}(K_1p-Lp'),\quad N=(-1)^h (K_1q-Lq').
\end{equation}
Thus the discriminant
$t^2-4n=T^2-4N$ is a quadratic expression in $K_1$ and~$L$.  

Note, however,
that many choices of $K_1$ and $L$ lead to inadmissible cases. For
instance, we must have $t^2-4n$ a positive non-square; and $Q$ positive
obviously constrains $L$. 

Additional data, such as a relationship on the
quantities
$P$, $Q$, and $Q'$, leads to just a family of 
polynomials quadratic in one variable as in the classical case; in that
context see Halter-Koch and Pacher
\cite{HP}, which expands considerably on \cite{118} and
\cite{111}.

The special case $Q=Q'$ entails $(-1)^{h+1}Q=K_1p'-\smash{q'}^2$. The
special case $P=P':=P_h$ is $Q\Div (2P+t)$; note that
$2P+t=(-1)^h(2qq'-Lp'-K_1p)$. Again in these cases,
the families of discriminants are quadratic \poly s in the surviving
parameter. Conversely, of course, the quadratic \poly s $D(L)$, say,
have the property that the period length of $\sqrt{D(l)}$, $l$
an integer for which $D(l)$ is positive and not a square, is constant
for all but possibly several exceptional~$l$. This case is discussed by
Schinzel~\cite{Schi}; see also
\cite{121} and note the `sleepers' of \cite{171}.
\begin{Example} \label{ex:4} Suppose we are given just $a_1=1$,
$a_2=4$; that is: $p=5$, $q=4$, $p'=1$, $q'=1$. Thus
\begin{alignat*}{3} P^2+tP+n &= Q(K_1-1) &\qquad &   A+P
&&= 4-5K_1\,,\\ A^2-tA+n &= 4K_2-Q &\qquad& 
A+P &&= Q-5K_2\,.
\end{alignat*} 
The equations on the right yield $Q-4=5(K_2-K_1)$. Subtracting
the first equation on the left from the second gives
$(A+P)(A-t-P)=4K_2-QK_1$. Substituting for $Q$ this is 
$(K_2-K_1)(4-5K_1)$, so indeed $A-t-P=K_2-K_1$. Thus if $Q=4$ then
$A-t=P$ and $2A-t=4-5K_1$. The smallest example then has $A-t=7$ with
$K_2=-2$ and so $t^2-4n=4\cdot61$.
\end{Example}
\begin{Example} \label{ex:5} Given $a_1=1$,
$a_2=4$, $a_3=3$, $a_4=1$, $a_5=2$, together with $Q_5=Q_6$ ---
so that we are half way in a period of odd length, we have
$p=58$, $q=47$, $p'=21$, $q'=17$. Thus 
\begin{alignat*}{3} -Q' &= 289-21K_1 &\qquad &   A+P
&&= 58K_1-799\,,\\ A^2-tA+n &= 289Q-47K_2 &\qquad& 
A+P &&= 58K_2-357Q\,.
\end{alignat*} 
On the right we see that $357Q-799\equiv9Q-45\equiv0\pmod{58}$
and only the smallest example, that with $Q=5$, is going to be
at all small. Moreover, if $Q=5$ then $K_2-K_1=q'K=17$, so $k=1$.
Recalling that
$Q'=Q$, the equation for
$-Q'$ gives $K_1=14$, and so $K_2=31$. All that amounts to $A+P=13$.
As above, we also have $(A+P)(A-t-P)$, so $2P+t=12$ and $2A-t=14$. We see
that  $t^2-4n=4\cdot61$. 

This example is just of the right size reliably to check our formulas just
above at the start of
\S\ref{ss:general}; we need only that the choice $Q=5$ is $L=38$, and then
$Q=Q'$ is $K_1=14$.

\end{Example}
\begin{Example} \label{ex:6}
Given $a_1=2$ and $P=P_2=P_3$, signalling halfway in a period
of even length, we have $p=2$, $q=1$, $p'=1$, $q'=0$. Hence $K_2=K_1$,
and
$(A+P)(A-t-P)=K_1(Q-1)$. Since $A+P=2K_1$ we see that $A-t-P=(Q-1)/2$ and $Q$
is odd. Thus the condition $Q\Div (2P+t)$ is $4K_1\equiv-1\pmod Q$. We obtain
a family of possible discriminants controlled by one near arbitrary
integer parameter and by $(2P+t)/Q$.
\end{Example}
\section{Function Fields} \label{s:4}
\noindent We will apply Theorem~\ref{th:2} in function fields (with base
field not of characteristic two) to determine \poly s of even degree whose
square root has a periodic
\cf\
\ex\ defined over the base field. We first look at quartic \poly s
$D(X)=A^2+4v(X+w)$ and, as example, first consider the case 
$$\sqrt D=\CF{A\\\overline{B\\C\\E\\C\\B\\2A}\,}$$
with $B$, $C$, and $E$ of degree~$1$, and $A$ of course of degree~$2$. Note
here that the
\emph{regulator}, namely the degree of the unit, or, equivalently, the sum of
the degrees of the
\pq s comprising the period, is $m=7$.

In order to apply the classical result we must first compute
\begin{multline*}
\begin{pmatrix} B&1\\1&0\end{pmatrix}\begin{pmatrix}
C&1\\1&0\end{pmatrix}
\begin{pmatrix} E&1\\1&0\end{pmatrix}
\begin{pmatrix} C&1\\1&0\end{pmatrix}\begin{pmatrix} B&1\\1&0\end{pmatrix}
\\=\begin{pmatrix} BCECB+ECB+2BCB+BCE+2B+E&BCEC+EC+2BC+1
\\CECB+2CB+CE+1&CEC+2C\end{pmatrix}.
\end{multline*}
Next we write say $B=b(X+\beta)$, $C=c(X+\gamma)$, and $E=e(X+\varepsilon)$,
allowing us explicitly to detail the \poly s $p$, $q=p'$, and $q'$ appearing in
the proof of Theorem~\ref{th:2}. In particular, we notice that $\deg p=5$,
$\deg q=4$, and $\deg q'=3$. We then endeavour to find $A$.

Indeed, by  Theorem~\ref{th:2} we know that $2A=Kp-qq'$ (noting that here
$t=0$ and $h=5$); in principle $K$ denotes an arbitrary \poly.

However, we require $\deg D=4$, and therefore that $\deg A=2$. Because 
$\deg qq'=7$ and $\deg p=5$ it follows that $K$ must be of degree $2$.
Thus we detail $Kp-qq'$, say with $K=kX^2+k'X+k''$, and set the 
leading coefficients of $2A$, those of $X^7$, $X^6$, $\ldots\,$, and $X^3$, equal
to zero.

Note that there are a total of $9$ unknowns $k$, $\ldots\,$, $b$, $\beta$,
$\ldots\,$, thus far controlled by
the $5$ conditions ensuring that $A$ is of degree just~$2$.

But we can tame $D(X)$, and hence $A(X)$, considerably further. First, it
loses no generality to divide by the leading coefficient of $D$
(this must be a square because $A$ is defined over the base field) and, thus,
to suppose
$D$ is monic. Second, we may translate
$X$ by a constant, thus fixing the coefficient of $X^3$ of $D$ (and hence the
coefficient of $X$ of
$A$), say to $0$. That leaves two degrees of freedom. Third, we may dilate
$X$, replacing it by a nonzero multiple, and then divide by the new leading
coefficient; one finds that for quartics such a dilation reduces the number
of free variables by one in all but the cases $m=2$ and $m=3$. From
related work \cite{158} it happens I know that a congenial dilation
is that yielding
$c=2$.

All that reduces the degrees of freedom to just one. We therefore should find
at most a one parameter family of possibilities for $A(X)=X^2+u$, say, thus for
$D$, and hence for each of the already mentioned `unknowns'.
    
\subsection{The Classical Result in Action}It would be neat to apply
the classical result to the example but, frankly, the hard yakka\footnote{
\textsf{yakka}: work [Australian Aboriginal].} involved seems inappropriate
here, so we'll retreat to the case $m=5$ and leave it to the reader to check
that there
$p=BCB+2B$, $q=p'=CB+1$, and $q'=B$, and that \emph{mutatis mutandis} we will
there too find a one parameter family of \cf s; in the case $m=5$ we will need
$\deg K=0$ and will set $K=k$.

For $m=5$, when we review the coefficients of $2A(X)=2(X^2+u)=Kp-qq'$ we find
that
$$\displaylines{X^3:kb^2c-bc^2 \quad\text{which vanishes, so
$kb=c$;}\quad\hfill\\
X^2:kb^2c(2\beta+\gamma)-b^2c(\beta+2\gamma)=bc^2(\beta-\gamma)\hfill\\
X^{\phantom{2}}:
kb^2c(2\beta\gamma+\beta^2)+2kb-bc^2(2\beta\gamma+\gamma^2)-c
=bc^2(\beta-\gamma)(\beta+\gamma)+c\hfill\\
\rlap{1}\hphantom{X^0}:
kb^2c\beta^2\gamma+2kb\beta-bc^2\beta\gamma^2-c\gamma
=bc^2(\beta-\gamma)\beta\gamma+c(2\beta-\gamma)
\hfill
}$$
Just so, it will be convenient to check the consistency of our solution by noting that the
coefficients of $Kq-{q^{\smash{\prime}}}^2=D-A^2=R$, say, are given by
$$\displaylines{X^2:kbc-c^2 \quad\text{which must vanish, so $kb=c$};\quad\hfill\\
X^{\phantom{2}}:kbc(\beta+\gamma)-2c^2\gamma=c^2(\beta-\gamma)\hfill\\
\rlap{1}\hphantom{X^0}:kbc\beta\gamma+k-c^2\gamma^2=c^2(\beta-\gamma)\gamma+k\hfill
}$$
Our normalisation $A=X^2+u$ gives $bc^2(\beta-\gamma)=2$, and then $2(\beta+\gamma)+c=0$.
Dilating
$X$ by the ratio $-(\beta+\gamma)$ is equivalent to taking $\beta+\gamma=-1$. If, further, we
choose to write $\beta-\gamma=s$, then $\beta=(s-1)/2$, $\gamma=-(s+1)/2$, and we obtain
$2u=(1-s^2)/2+(3s-1)$, or $u=-(s^2-6s+1)/4$. Because $bc^2(\beta-\gamma)=4bs=2$  we have
$b=1/2s$, $c=2$, and $k=4s$.

Further, we find that $R(X)=4sX-2s(s+1)+4s=4sX-2s(s-1)$. Thus
\begin{equation}\label{eq:fam5}
Y^2=D(X)=\bigl(X^2-\tfrac14(s^2-6s+1)\bigr)^2+4s\bigl(X-\tfrac12(s-1)\bigr)
\end{equation}
is the family of monic quartic \poly s defined over $\K=\Q(s)$ and with zero
trace so that the function field $\K(X,Y)$ has a unit of regulator~$5$.

One should promptly check such an allegation. Indeed, we find that 
\begin{multline*}
Y=\lCF{X^2-\tfrac14(s^2-6s+1)\\\ol{\bigl(X+\tfrac12(s-1)\bigr)/2s}}\\
\rCF{\ol{2\bigl(X-\tfrac12(s+1)\bigr)\\\bigl(X+\tfrac12(s-1)\bigr)/2s
\\2\bigl(X^2-\tfrac14(s^2-6s+1)\bigr)}}.
\end{multline*}
This computation provides the reminder that of course $s=0$ is not admissible because, if $s=0$,
then
$D$ is a square.

\subsection{The Generalised Result in Action} \label{ss:PlanA'}
Notwithstanding this success of the classical result, the case
$m=7$ already seems too painful for any other than a willing and 
energetic student. Our generalisation should however at least halve the
pain and, it is to be hoped, more than double the willingness. 

Recall we suppose $Y^2=D(X)=(X^2+u)^2+4v(X+w)=A^2+R$. In applying the new
result we first compute just
$$\begin{pmatrix} B&1\\1&0\end{pmatrix}\begin{pmatrix} C&1\\1&0\end{pmatrix}
=\begin{pmatrix} BC+1&B\\C&1\end{pmatrix}=\begin{pmatrix} p&p'\\q&q'\end{pmatrix}
$$
so $p=BC+1$, $p'=B$, $q=C$, $q'=1$ and note that $Q\Div(2P+t)$. 
Checking degrees in the equations at \S\S\ref{ss:general} we see the $K_i$
are constants
$k_i$, and as remarked, we obtain
$A-t-P=k_2-k_1=k$, say, $(BC+1)k=(BQ-C)$, $A+P=C-k_1(BC+1)$. 

Comparing coefficients in that last equation yields
$$\displaylines{
X^2: 2=-k_1bc\hfill\\
X\phantom{^2}: 0=
c-k_1bc(\beta+\gamma)=c+2(\beta+\gamma),\hfill }$$
the latter by normalisation whereby $A$, and so also $P$, has zero trace. By
dilation we may choose $c=2$ obtaining $k_1b=-1$ and $\beta+\gamma=-1$.

We note that $E$ is the integer part of $(Y+P)/Q$. Hence, because $D(X)$
is monic and has  zero trace, necessarily $Q=2(X-\varepsilon)/e$.
Thus
$(BC+1)k=(BQ-C)$ alleges that
$$\displaylines{
X^2: kbc=2b/e, \quad\text{so}\quad 2=kce;\hfill\\
X\phantom{^2}: kbc(\beta+\gamma)=2b(\beta-\varepsilon)/e-c
\quad\text{or}\quad kb(\gamma+\varepsilon)=-1;\hfill\\
\rlap{1}\hphantom{X^0}: kbc\beta\gamma+k=-2b\beta\varepsilon/e-c\gamma
\quad\text{or}\quad kbc\beta(\gamma+\varepsilon)=-k-c\gamma=-c\beta\,,\hfill
\\\hfill\quad\text{and so}\quad c(\beta-\gamma)=k
.\quad}$$
The equation $-Q'=(-1)^{h+1}(q'{}^2-K_1p')$ contains no new information
other than for the reminder that in fact $C=c(X+\gamma)$ entails
$Q'=2(X-\gamma)/c$. Just so the equation for $A^2-tA+n$ simply gives 
$-R$.

However, the condition that $Q=2(X-\varepsilon)/e$ divides
$2P+t=-k+C-k_1BC-k_1$ says that
$-k+c(X+\gamma)-k_1bc(X^2+(\beta+\gamma)X+\beta\gamma)-k_1$
has $\varepsilon$ as a zero.
That is
$$-k_1bc\varepsilon^2=k-c\gamma+k_1bc\beta\gamma+k_1 \quad\text{so}\quad
2\varepsilon^2=k-(2\gamma+2\beta\gamma-k_1).$$
Subtracting $2\gamma^2$ from each side of this equation we get
$$2(\varepsilon^2-\gamma^2)=2(\varepsilon-\gamma)(\varepsilon+\gamma)
=k-2\gamma-2\gamma(\beta+\gamma)+k_1
=k+k_1\,.$$
Set $\varepsilon+\gamma=-s$. Notice that $kbs=1$ and $k_1b=-1$ gives
$k_1=-ks$. We have
$$
\begin{aligned}
\beta+\gamma&=-1\qquad\text{or}\qquad&&2(\beta+\gamma)s=-2s\\
\varepsilon+\gamma&=-s\qquad&&2(\varepsilon+\gamma)s=-2s^2\\
2(\beta-\gamma)&=k\qquad&&2(\beta-\gamma)s=ks\\
2(\varepsilon-\gamma)s&=k(s-1)\qquad&&2(\varepsilon-\gamma)s=k(s-1)\,.\\
\end{aligned} 
$$
It is now straightforward to solve the various equations and to obtain
\begin{align}
B&=\bigl(X+\tfrac12(s^2-s-1)\bigr)/2s^2(s-1)\notag\\\label{eq:solution}
C&=2\bigl(X-\tfrac12(s^2-s+1)\bigr)\\
E&=\bigl(X+\tfrac12(s^2-3s+1)\bigr)/2s(s-1)\notag\,,
\end{align}
as well as $k=2s(s-1)$.   As remarked, we obtain $A$ from
$2A-t=k+C-k_1BC-k_1$ and $v(X+w)$ as $Q-(k+k_1)C$. Thus
$Y^2=(X^2+u)^2+4v(X+w)=A^2+R$ with
$$A=X^2-\tfrac14(s^4-6s^3+3s^2+2s+1) \quad\text{and}\quad
R=4s^2(s-1)\bigl(X-\tfrac12(s^2-s-1)\bigr).
$$
On our checking, a direct computation of the \cfe\ of $Y$ so given indeed
yields the 
\pq s predicted at \eqref{eq:solution} above.

\section{Comments}

\subsection{Schinzel's theorem} Andrzej Schinzel~\cite{Schi} shows that
if a \poly\ $D(X)$ taking integer values at integers $l$ has the
property that the length of the period of the \cfe\ of $\sqrt{D(l)}$ is
bounded as $l\to\infty$ then (a) the function field $\Q(X,Y)$ ---
where $Y=\sqrt{D(X)}$ --- contains non-trivial units; equivalently, the
function field \cfe\ of $Y$ is periodic, and (b) some nontrivial units
$a(X)+b(X)Y$ in $\Q(X,Y)$ have both $a$ and $b$ in $\Z[X]$. If $\deg
D=2$, say $D(X)=A^2X^2+BX+C$, then the function field condition (a) is
trivial, but the arithmetic condition (b) entails that the
discriminant $(B^2-4A^2C)$ divides
$4 (2A^2,B)^2$.

Section~\ref{s:4} above provides examples of classes of \poly s of
degree~$4$ satisfying the function field condition (a). Roger Patterson,
see~\cite{Pa}, has carried through the interesting exercise of finding
which values of the parameters lead to (b) also being satisfied; see also
comments in~\cite{171}.

\subsection{Explicit continued fraction expansion} I have learned how
explicitly to expand the general quartic (and sextic)  \poly, that is,
I found nontrivial recursion relations on the complete
quotients; see~\cite{163}. In principle, at any rate, those
techniques generalise to higher genus cases. Thus, fortunately perhaps,
the ideas of
\S\ref{s:4} are mostly an amusing diversion rather than a necessary
method.

\subsection{Short periods and long periods} In principle one is
interested in relatively long periods, in the hope of obtaining insight
into Gauss's conjecture to the effect that the majority of real quadratic
number fields with prime discriminant have class number one. In that
context, short periods, as discussed here only give one a view of the
enemy. Nonetheless, as wonderfully exemplified by recent work of
Andr\'as Bir\'o, \cite{Bi1, Bi2}, the short case raises
fascinating issues. The \poly s
$X^2+4$ and
$4X^2+1$ satisfy Schinzel's conditions so we can write the corresponding
numerical periods explicitly. That gives, as Bir\'o ingeniously shows,
sufficient handle on the matter after all to determine all cases $X=l$
that yield class number one. But, for example, the corresponding class
number two problems appear yet less accessible because as $l$ varies one
seems not to have a bound on the period length of expansions of 
reduced elements not from the principal class.

It's not so much the length of the period as the size of the
regulator (the logarithm of the absolute value of the fundamental
unit) that matters. Discriminants $D$ belonging to families alluded
to in this paper have regulator $O(\log D)$ only. In fact, see
\cite{Pa, 171}, it seems that if one has a formula for the units of a
parametrised family then the regulators are never more than 
$O\bigl((\log D)^2\bigr)$.


\bibliographystyle{amsalpha}

\label{page:lastpage}
\end{document}